\documentclass[11pt, twoside]{article}
\usepackage{amssymb}
\usepackage{amsmath, amssymb}\usepackage{cite}\usepackage{mathrsfs}\usepackage[amsmath, thmmarks]{ntheorem}
\usepackage{listings}
\usepackage[titletoc]{appendix}\allowdisplaybreaks
\usepackage{multirow}

\makeatletter
\renewcommand\theequation{\thesection.\@arabic\c@equation}
\makeatother

\newtheorem{thm}{Theorem}[section]%
\newtheorem{lem}[thm]{Lemma}%
\newtheorem{Con}[thm]{Conjecture}%
\input cyracc.def

\def\f{\noindent}

\def\demo{\f{\bf Proof}\hskip10pt}

\def\qed{\hfill $\Box$}

\begin{document}

\title{{\bf On the permutations that strongly avoid the pattern 312 or 231}}
\footnotetext{The work was supported by the National Natural Science Foundation of China (No. 12061030) and Hainan Provincial Natural Science Foundation of China (No. 122RC652).\\
$^{*}$ Corresponding author.\\E-mail addresses: Junyao$_{-}$Pan@126.com (Junyao Pan); guopf999@163.com (Pengfei Guo).}

\author{{\bf Junyao Pan$^{{\rm 1}}$, Pengfei Guo$^{{\rm 2}, *}$}\\
{\footnotesize 1. School of Sciences, Wuxi University, Wuxi, Jiangsu, 214105, P. R. China}\\
{\footnotesize 2. School of Mathematics and Statistics, Hainan Normal University, Haikou, Hainan 571158, P. R. China}}

\date{}
\maketitle

%
%
%
%

\noindent{\small {\bf Abstract:} In 2019, B\'ona and Smith introduced the notion of \emph{strong pattern avoidance},
that is, a permutation and its square both avoid a given pattern. In this paper, we enumerate the set of permutations $\pi$
which not only strongly avoid the pattern $312$ or $231$ but also avoid the pattern $\tau$, for $\tau\in S_3$ and some $\tau\in S_4$.
One of them is to give a positive answer to a conjecture of Archer and Geary.

\vskip0.2cm
\noindent{\small {\bf Keywords:} Permutation; strong pattern avoidance; pattern $312$; pattern $231$.


\section {Introduction}

Let $S_n$ denote the symmetry group on $[n]:=\{1, 2, \ldots , n\}$. It is well-known that every permutation $\pi$ in $S_n$
can be written in its one-line notation as $\pi = \pi_1\pi_2\cdots \pi_n$, where $\pi_i = \pi(i)$ for all $i \in [n]$.
From the perspective, the patterns can be considered in the permutation. Let $\pi= \pi_1\pi_2\cdots \pi_n\in S_n$
and $\tau= \tau_1\tau_2\cdots \tau_k\in S_k$ with $k\leq n$. If there exists a subset of indices $1\leq i_1<i_2<\cdot\cdot\cdot<i_k\leq n$
such that $\pi_{i_s}>\pi_{i_t}$ if and only if $\tau_s>\tau_t$ for all $1\leq s<t\leq k$, then we say that $\tau$ is \emph{contained} in $\pi$
and the subsequence $\pi_{i_1}\pi_{i_2}\cdot\cdot\cdot \pi_{i_k}$ is called an \emph{occurrence} of $\tau$ in $\pi$ and denoted by $\tau\leq \pi$.
For example, $132\leq24153$, because $2,~5,~3$ appear in the same order of size as the letters in $132$.
The permutation that does not contain an occurrence of $\tau$ is called a $\tau$-\emph{avoiding} permutation.
The theory of pattern avoidance in permutations was introduced by Knuth in \cite{K}, which has been widely studied for half a century, refer to \cite{B, V}.
Recently, the research of the pattern avoidance concerning the symmetric group concept has attracted the attention of some scholars.
However, it becomes more complicated as the conditions strengthen and the results about this topic can be seen in \cite{AE, AG, BC, BS, BD, K, P, SS}.

In 2019, B\'ona and Smith proposed the definition of \emph{strong pattern avoidance} in \cite{BS},
that is, a permutation $\pi$ strongly avoids $\tau$ if both $\pi$ and $\pi^2$ avoid $\tau$. Shortly afterward,
Burcroff and Defant \cite{BD} introduced the definition of \emph{powerful avoidance}, namely,
a permutation $\pi$ powerfully avoids $\tau$ if every power of $\pi$ avoids $\tau$. More recently, Archer and Geary \cite{AG}
generalized the ideas of strongly avoidance and powerful avoidance to what they called \emph{chain avoidance},
that is, a permutation $\pi$ avoids a chain of patterns $(\tau_1 : \tau_2 : \cdot\cdot\cdot : \tau_k)$
if the $i$-th power of $\pi$ avoids the pattern $\tau_i$. In addition, this definition can be extended to the set of patterns as well,
which is separated by a comma. For example, if $\pi$ avoids the chain $(\sigma_1, \sigma_2 : \tau)$, then $\pi$ avoids both $\sigma_1$ and $\sigma_2$,
and $\pi^2$ avoids $\tau$. Let $C_n(\sigma_1, \sigma_2 : \tau)$ denote the set of permutations in $S_n$ that avoid the chain $(\sigma_1, \sigma_2 : \tau)$,
and $c_n(\sigma_1, \sigma_2 : \tau)$ denote the number of such permutations. It is worth mentioning that
Archer and Geary \cite{AG} not only gave a formula of computing the number of
unimodal permutations whose square avoids every $\tau\in S_3$ but also proposed an interesting conjecture as shown below:

\begin{Con}\label{pan1-1}\normalfont([2])
 $c_n(231,1432 : 231)=L_{n+1}-\lceil \frac{n}{2}\rceil-1$ for all positive integer $n$, where $L_{n+1}$ is the $(n+1)$-th Lucas number.
\end{Con}

In this paper, we not only give a positive answer to the Conjecture\ \ref{pan1-1} but also enumerate $C_n(231, \tau : 231)$ and $C_n(312, \tau : 312)$ for $\tau\in S_3$ and some $\tau\in S_4$. A summary of these results can be found in the following table.

\begin{center}
\begin{tabular}{|c|c|c|c|}
\hline
$\tau_{1}$ &$\tau_{2}$ & $c_n(231,\tau_{1}: 231)=c_n(312,\tau_{2}: 312)$ & Theorem  \\
\hline\hline
123 & 123 & $2n-3$ & Theorem \ref{pan3-1} \\
\hline
321 & 321 & $F_{n}$ & Theorem \ref{pan3-2} \\
\hline
132 & 213 & $k^2+1$ if $n =2k$; $k^2+k+1$ if $n =2k+1$ & Theorem \ref{pan3-3} \\
\hline
312 & 231 & $2^{n-1}$ & Theorem \ref{pan3-4} \\
\hline
213 & 132 & $k^2+1$ if $n =2k$; $k^2+k+1$ if $n=2k+1$ & Theorem \ref{pan3-5} \\
\hline
1432 & 3214 & $L_{n+1}-\lceil \frac{n}{2}\rceil-1$ & Theorem \ref{pan4-1} \\
\hline
1423 & 2314 & $\frac{7}{3}\cdot4^{k-1}-\frac{1}{3}$ if $n=2k$; $\frac{14}{3}\cdot4^{k-1}-\frac{2}{3}$ if $n=2k+1$ & Theorem \ref{pan4-2} \\
\hline
1243 & 2134 & $\frac{4k^3+3k^2-k}{6}+1$ if $n=2k$; $\frac{4k^3+9k^2+5k}{6}+1$ if $n=2k+1$ & Theorem \ref{pan4-3} \\
\hline
2143 & 2143 & $\frac{4k^3+3k^2-k}{6}+1$ if $n =2k$; $\frac{4k^3+9k^2+5k}{6}+1$ if $n =2k+1$ & Theorem \ref{pan4-4} \\
\hline
\end{tabular}
\end{center}

\section {Preliminaries}
Let $\pi = \pi_1\pi_2\cdots \pi_n\in S_n$. Recall that the reverse of $\pi$ is the permutation $rev(\pi)=\pi_n\pi_{n-1}\cdots \pi_1$,
the complement of $\pi$ is the permutation $comp(\pi)=(n+1-\pi_1)(n+1-\pi_2)\cdots (n+1-\pi_n)$, and the reverse complement of $\pi$ is the permutation $comp(rev(\pi))$.
In addition, let $S$ be a set of permutations. We use $comp(rev(S))$ to denote the set $\{comp(rev(\pi))|\pi\in S\}$.
Inspired by B\'ona and Smith, Burcroff and Defant observed that $comp(rev(\pi^k))=comp(rev(\pi))^k$ for any positive integer $k$.
Therefore, it is easy to deduce the following lemma.

\begin{lem}\label{pan2-1}\normalfont
Let $(S_1 : S_2 : \cdot\cdot\cdot : S_k)$ and $(S'_1 : S'_2 : \cdot\cdot\cdot : S'_k)$ be two chains of patterns
such that $S'_i=comp(rev(S_i))$ for $i=1,2, \ldots , k$. Then $c_n(\tau_1 : \tau_2 : \cdot\cdot\cdot : \tau_k)=c_n(\tau'_1 : \tau'_2 : \cdot\cdot\cdot : \tau'_k)$.
\end{lem}
Actually, applying Lemma \ref{pan2-1}, we see that all results of Archer and Geary \cite{AG} can be attributed to their inverse complements,
such as $c_n(132,231:123)=\lfloor\frac{n}{2}\rfloor$ for $n\geq8$, $c_n(132,231:213)=n+1$ for $n\geq3$, and so on. Note that $312=comp(rev(231))$ and $3214=comp(rev(1432))$.
Therefore, $c_n(231,1432 : 231)=c_n(312,3214 : 312)$. Additionally, B\'ona and Smith \cite{BS} obtained
a result about strongly avoidance of the pattern $312$ that will be repeatedly applied, as follows:

\begin{thm}\label{pan2-2}\normalfont([5, Theorem 3.1])
For any permutation $\pi$ ending in $1$, the following two statements are equivalent.

(A) The permutation $\pi$ is strongly $312$-avoiding.

(B) The permutation $\pi$ has form $\pi=(k+1)(k+2)\cdot\cdot\cdot nk(k-1)(k-2)\cdot\cdot\cdot1$ where $k\geq\frac{n}{2}$.

That is, $\pi$ is unimodal beginning with its $n-k\leq\frac{n}{2}$ largest entries in increasing order followed

by the remaining $k$ smallest entries in decreasing order.
\end{thm}
Thereby, we enumerate the set of permutations $\pi$ which not only strongly avoid the pattern $312$ but also avoid the pattern $\tau$, for $\tau\in S_3$ and some $\tau\in S_4$. In particular, we shall give a positive answer to the Conjecture\ \ref{pan1-1} by computing $c_n(312,3214 : 312)$.

Next we give a simple while important formula.

\begin{lem}\label{pan2-3}\normalfont
If $n$ is a positive integer, then
\begin{align*}
\sum_{i=1}^{n}\lceil \frac{i}{2}\rceil&=
  \begin{cases}
    k^2+k & \text{if } n =2k, \\
    k^2+2k+1 & \text{if } n =2k+1.
  \end{cases}
\end{align*}
\end{lem}
\demo We claim that $\lceil \frac{n-1}{2}\rceil+\lceil \frac{n}{2}\rceil=n$. If $n=2k$, we see that $\lceil \frac{n-1}{2}\rceil=k$ and $\lceil \frac{n}{2}\rceil=k$, and thus $\lceil \frac{n-1}{2}\rceil+\lceil \frac{n}{2}\rceil=2k=n$; and if $n=2k+1$, then we have $\lceil \frac{n-1}{2}\rceil=k$ and $\lceil\frac{n}{2}\rceil=k+1$, and so $\lceil \frac{n-1}{2}\rceil+\lceil \frac{n}{2}\rceil=2k+1=n$. Therefore, our claim is true. Now we start to prove this lemma. If $n=2k$, then $\sum_{i=1}^{n}\lceil \frac{i}{2}\rceil=(\lceil\frac{1}{2}\rceil+\lceil\frac{2}{2}\rceil)+\cdot\cdot\cdot+(\lceil\frac{n-1}{2}\rceil+\lceil\frac{n}{2}\rceil)$. By using our claim, we deduce that $$\sum_{i=1}^{n}\lceil \frac{i}{2}\rceil=2+4+\cdot\cdot\cdot+2k=k^2+k.$$ If $n=2k+1$, then $\sum_{i=1}^{n}\lceil \frac{i}{2}\rceil=\sum_{i=1}^{2k}\lceil \frac{i}{2}\rceil+\lceil \frac{2k+1}{2}\rceil$ and thus $\sum_{i=1}^{n}\lceil \frac{i}{2}\rceil=k^2+2k+1$. The proof of this lemma is completed.  \qed

Finally, we introduce some notations that will be used. Given $\sigma\in S_l$ and $\mu\in S_m$. Let $\sigma\oplus\mu$ and $\sigma\ominus\mu$ denote the direct sum and skew sum of $\sigma$ and $\mu$ respectively, those are,
\begin{align*}
\sigma\oplus\mu(i)&=
  \begin{cases}
    \sigma(i) & \text{if } 1\leq i\leq l; \\
    \mu(i)+l & \text{if } l+1\leq i\leq l+m
  \end{cases}
\end{align*}
and
\begin{align*}
\sigma\ominus\mu(i)&=
  \begin{cases}
    \sigma(i)+m & \text{if } 1\leq i\leq l; \\
    \mu(i-l) & \text{if } l+1\leq i\leq l+m.
  \end{cases}
\end{align*}
Let $C^{(i)}_n(S_1 : S_2 : \cdot\cdot\cdot : S_k)=\{\pi\in C_n(S_1 : S_2 : \cdot\cdot\cdot : S_k)|\pi(i)=1\}$ and $c^{(i)}_n(S_1 : S_2 : \cdot\cdot\cdot : S_k)=|C^{(i)}_n(S_1 : S_2 : \cdot\cdot\cdot : S_k)|$ for every $i\in[n]$, where $(S_1 : S_2 : \cdot\cdot\cdot : S_k)$ is a chain of patterns. Note that $$c_n(S_1 : S_2 : \cdot\cdot\cdot : S_k)=\sum_{i=1}^{n}c^{(i)}_n(S_1 : S_2 : \cdot\cdot\cdot : S_k).$$ Indeed, this idea will be an effective means for us to provide calculating formulas. Moreover, let $\delta_n=n\cdot\cdot\cdot321$ be the reverse of the identity $id_n=123\cdot\cdot\cdot n$.

\section {Avoiding pattern of length 3}

In this section, we investigate $c_n(231,\tau : 231)$ for $\tau\in S_3$.

\begin{thm}\label{pan3-1}\normalfont
If $n\geq3$, then $$c_n(231,123 : 231)=c_n(312,123 : 312)=2n-3.$$
\end{thm}
\demo Since $312=comp(rev(231))$ and $123=comp(rev(123))$, it follows from Lemma \ref{pan2-1} that $c_n(231,123 : 231)=c_n(312,123 : 312)$. Consider $c_n(312,123 : 312)$. It is straightforward to see that $c_3(312,123 : 312)=3$. Next we consider $n>3$. Note that $$c_n(312,123 : 312)=\sum_{i=1}^{n}c^{(i)}_n(312,123 : 312).$$
If $\pi\in C^{(n)}_n(312,123 : 312)$, then we see that either $\pi=n(n-1)\cdots 21$ or $\pi=(n-1)n(n-2)\cdots 21$ by Theorem\ \ref{pan2-2}, and thus $c^{(n)}_n(312,123 : 312)=2$. If $\pi\in C^{(1)}_n(312,123 : 312)$, then $\pi=1n(n-1)\cdots 2$ because $\pi$ avoids $123$. In this case, $\pi^2=123\cdots n$, and so $c^{(1)}_n(312,123 : 312)=1$. Similarly, we infer that $\pi\in C^{(2)}_n(312,123 : 312)$ if and only if $\pi=21n(n-1)\cdots 3$, and hence $c^{(2)}_n(312,123 : 312)=1$.
Let $\pi = \pi_1\pi_2\cdots \pi_n\in S_n$ with $\pi_i=1$ and $2<i\leq n-1$. Clearly, $\pi\in C^{(i)}_n(312,123 : 312)$ if and only if $\pi=\sigma\oplus\delta_{n-i}$ where $\sigma\in C^{(i)}_i(312,123 : 312)$. Moreover, $c^{(3)}_3(312,123 : 312)=1$ and $c^{(i)}_i(312,123 : 312)=2$ for $3<i\leq n-1$. Thus, $c^{(3)}_n(312,123 : 312)=1$ and $c^{(i)}_n(312,123 : 312)=2$ for $3<i\leq n-1$.

According to the above discussions, we deduce that $c_n(312,123 : 312)=2n-3$, as desired. \qed

\begin{thm}\label{pan3-2}\normalfont
Let $F_{n}$ be the $n$-th Fibonacci number. If $n\geq1$, then $$c_n(231,321 : 231)=c_n(312,321 : 312)=F_n.$$
\end{thm}
\demo By the same token, we deduce that $c_n(231,321 : 231)=c_n(312,321 : 312)$. So it suffices to consider $c_n(312,321 : 312)$. It is easy to check that $c_1(312,321 : 312)=1$, $c_2(312,321 : 312)=2$ and $c_3(312,321 : 312)=3$. Next we discuss $c_n(312,321 : 312)$ for $n>3$.

By Theorem\ \ref{pan2-2}, we note that $c^{(n)}_n(312,321:312)=0$. Additionally, it is obvious that $\pi \in C^{(1)}_n(312,321 : 312)$ if and only if $\pi=1\oplus\mu$ where $\mu\in C_{n-1}(312,321 : 312)$, and therefore $c^{(1)}_n(312,321 : 312)=c_{n-1}(312,321 : 312)$. Similarly, we deduce that $c^{(2)}_n(312,321 : 312)=c_{n-2}(312,321 : 312)$. Let $\pi = \pi_1\pi_2\cdots \pi_n\in S_n$ with $\pi_i=1$ and $2<i\leq n-1$. If $\pi$ avoids $321$ and $312$, then $\pi=(id_{i-1}\ominus1)\oplus\mu$ where $\mu$ avoids $321$ and $312$. However, $\pi^2$ contains the pattern $312$, and thus $c^{(i)}_n(312,321 : 312)=0$ for $2<i\leq n-1$. So far, we have seen that $$c_n(231,321 : 231)=c_{n-1}(312,321 : 312)+c_{n-2}(312,321 : 312).$$
By the definition of Fibonacci number again, the result holds immediately.    \qed

\begin{thm}\label{pan3-3}\normalfont
If $n\geq1$, then
\begin{align*}
c_n(231,132 : 231)=c_n(312,213 : 312)&=
  \begin{cases}
    k^2+1 & \text{if } n =2k, \\
    k^2+k+1 & \text{if } n =2k+1.
  \end{cases}
\end{align*}
\end{thm}
\demo In a similar manner, we see that $c_n(231,132 : 231)=c_n(312,213 : 312)$. So it suffices to confirm $c_n(312,213 : 312)$. It is easy to check that $c_1(312,213 : 312)=1$, $c_2(312,213 : 312)=2$ and $c_3(312,213 : 312)=3$. Next we discuss $c_n(312,213 : 312)$ for $n>3$.

 By Theorem\ \ref{pan2-2}, we deduce that $c^{(n)}_n(312,213 : 312)=n-\lceil \frac{n}{2}\rceil$. In addition, it is clear that $\pi\in C^{(1)}_n(312,213 : 312)$ if and only if $\pi=1\oplus\mu$ where $\mu\in C_{n-1}(312,213 : 312)$, and thus $c^{(1)}_n(312,213 : 312)=c_{n-1}(312,213 : 312)$. Let $\pi = \pi_1\pi_2\cdots \pi_n\in S_n$ with $\pi_i=1$ and $1<i\leq n-1$. In this case, it is clear that $\pi$ either contains $312$ or contains $213$. Therefore, $c^{(i)}_n(312,213 : 312)=0$ for $1<i\leq n-1$. So far, we have seen that, $$c_n(312,213 : 312)=c_{n-1}(312,213 : 312)+n-\lceil \frac{n}{2}\rceil.$$
In fact, the above equation holds for all $n>1$. Thereby, we deduce that $$c_n(312,213 : 312)-c_1(312,213 : 312)=\sum_{i=2}^{n}i-\sum_{i=2}^{n}\lceil \frac{i}{2}\rceil.$$
According to Lemma\ \ref{pan2-3}, we deduce that
\begin{align*}
c_n(312,213 : 312)&=
  \begin{cases}
    k^2+1 & \text{if } n =2k, \\
    k^2+k+1 & \text{if } n =2k+1.
  \end{cases}
\end{align*}
The proof of this theorem is completed.    \qed

\begin{thm}\label{pan3-4}\normalfont
If $n\geq 1$, then $$c_n(231,312 : 231)=c_n(312,231 : 312)=2^{n-1}.$$
\end{thm}
\demo By Lemma\ \ref{pan2-1}, we infer that $c_n(231,312 : 231)=c_n(312,231 : 312)$. So it suffices to consider $c_n(312,231 : 312)$. One easily checks that $c_1(312,231 : 312)=1$, $c_2(312,231 : 312)=2$, $c_3(312,231 : 312)=4$ and $c_4(312,231 : 312)=8$. Next we discuss $c_n(312,231 : 312)$ for $n>3$.

It follows from Theorem\ \ref{pan2-2} that $c^{(n)}_n(312,231 : 312)=1$. Note that $\pi\in C^{(1)}_n(312,231 : 312)$ if and only if $\pi=1\oplus\mu$ where $\mu\in C_{n-1}(312,231 : 312)$. Thus, $c^{(1)}_n(312,231 : 312)=c_{n-1}(312,231 : 312)$. Let $\pi = \pi_1\pi_2\cdots \pi_n\in S_n$ with $\pi_i=1$ and $1<i\leq n-1$. Similarly, $\pi\in C^{(i)}_n(312,231 : 312)$
 if and only if $\pi=\delta_i\oplus\mu$ where $\mu\in C_{n-i}(312,231 : 312)$, and so $c^{(i)}_n(312,231 : 312)=c_{n-i}(312,231 : 312)$ for $1<i\leq n-1$. Therefore, we infer that for $n>3$, $$c_n(312,231 : 312)=\sum_{i=1}^{n-1}c_{i}(312,231 : 312)+1.$$
Indeed, the above equation holds for all $n\geq1$. Thereby, we deduce that $c_n(312,231 : 312)=2^{n-1}$, as desired.    \qed

\begin{thm}\label{pan3-5}\normalfont
If $n\geq1$, then
\begin{align*}
c_n(231,213 : 231)=c_n(312,132 : 312)&=
  \begin{cases}
    k^2+1 & \text{if } n =2k, \\
    k^2+k+1 & \text{if } n =2k+1.
  \end{cases}
\end{align*}

\end{thm}
\demo It follows from Lemma\ \ref{pan2-1} that $c_n(231,213 : 231)=c_n(312,132 : 312)$. Thus, it suffices to consider $c_n(312,132 : 312)$. One easily checks that $c_1(312,132 : 312)=1$, $c_2(312,132 : 312)=2$ and $c_3(312,132 : 312)=3$. Thus, the theorem is true for $n=1,2,3$, and next we it for $n>3$.

Applying Theorem\ \ref{pan2-2}, we see that $c^{(n)}_n(312,132 : 312)=n-\lceil \frac{n}{2}\rceil$. Furthermore, it is clear that $\pi\in C^{(1)}_n(312,132 : 312)$ if and only if $\pi=id_n$, and thus $c^{(1)}_n(312,132 : 312)=1$. Similarly, $\pi = \pi_1\pi_2\cdots \pi_n\in C^{(2)}_n(312,132 : 312)$ if and only if $\pi=213\cdots n$, and so $c^{(2)}_n(312,132 : 312)=1$. Let $\pi = \pi_1\pi_2\cdots \pi_n\in S_n$ with $\pi_i=1$ and $2<i\leq n-1$. Note that $\pi\in C^{(i)}_n(312,132 : 312)$ if and only if $\pi=\mu\oplus id_{n-i}$ where $\mu\in C^{(1)}_{i}(312,132 : 312)$, and so $c^{(i)}_n(312,132 : 312)=i-\lceil \frac{i}{2}\rceil$ for $2<i\leq n-1$. So far, we have seen that for $n>3$, $$c_n(312,132 : 312)=\sum_{i=3}^{n}(i-\lceil \frac{i}{2}\rceil)+2.$$
Therefore, for $n\geq4$, we deduce that
$$c_n(312,132 : 312)=(\sum_{i=2}^{n}i-2)-(\sum_{i=2}^{n}\lceil \frac{i}{2}\rceil-1)+2=\sum_{i=2}^{n}i-\sum_{i=2}^{n}\lceil \frac{i}{2}\rceil+1.$$
Note that $c_n(312,132 : 312)=c_n(312,213 : 312)$, as desired.    \qed

\section {Avoiding pattern of length 4}
In this section, we investigate $c_n(231,\tau : 231)$ for some $\tau\in S_4$. Next we start from giving a positive answer to Conjecture \ref{pan1-1}.

\begin{thm}\label{pan4-1}\normalfont
Let $L_{n+1}$ be the $(n+1)$-th Lucas number. If $n\geq1$, then $$c_n(231,1432 : 231)=c_n(312,3214 : 312)=L_{n+1}-\lceil \frac{n}{2}\rceil-1.$$
\end{thm}
\demo Since $312=comp(rev(231))$ and $3214=comp(rev(1432))$, then $c_n(231,1432 : 231)=c_n(312,3214 : 312)$ by Lemma \ref{pan2-1}. Thus, it suffices to verify $c_n(312,3214 : 312)=L_{n+1}-\lceil \frac{n}{2}\rceil-1$. Moreover, it is obvious that $$c_n(312,3214 : 312)=\sum_{i=1}^{n}c^{(i)}_n(312,3214 : 312).$$
Next we discuss $c^{(i)}_n(312,3214 : 312)$ in three different situations.

If $i=n$ then $c^{(n)}_n(312,3214 : 312)=n-\lceil \frac{n}{2}\rceil$ by Theorem\ \ref{pan2-2}. In addition, it is easy to see that $c^{(1)}_n(312,3214 : 312)=c_{n-1}(312,3214 : 312)$ and $c^{(2)}_n(312,3214 : 312)=c_{n-2}(312,3214 : 312)$. Let $\pi = \pi_1\pi_2\cdots \pi_n\in C^{(i)}_n(312,3214 : 312)$ for $2<i\leq n-1$. Since $\pi$ avoids $312$ and $3214$, we have that $\pi_1=2,\pi_2=3, \ldots ,\pi_{i-1}=i$. However, it is obvious that $\pi^2$ contains the pattern $312$, a contradiction. Therefore, $c^{(i)}_n(312,3214 : 312)=0$ for $2<i\leq n-1$. So far, we have found that $$c_n(312,3214 : 312)=c_{n-1}(312,3214 : 312)+c_{n-2}(312,3214 : 312)+n-\lceil \frac{n}{2}\rceil.$$

One easily checks that $c_n(312,3214 : 312)=L_{n+1}-\lceil \frac{n}{2}\rceil-1$ holds for $n=1, 2, 3$.
Proof by induction on $n$. Since $c_n(312,3214 : 312)=c_{n-1}(312,3214 : 312)+c_{n-2}(312,3214 : 312)+n-\lceil \frac{n}{2}\rceil$, we deduce that $$c_n(312,3214 : 312)=L_{n}-\lceil \frac{n-1}{2}\rceil-1+L_{n-1}-\lceil \frac{n-2}{2}\rceil-1+n-\lceil \frac{n}{2}\rceil.$$
Note that $\lceil \frac{n-1}{2}\rceil+\lceil \frac{n-2}{2}\rceil+1=n$. Therefore, $c_n(312,3214 : 312)=L_{n}+L_{n-1}-\lceil \frac{n}{2}\rceil-1$, as desired. \qed

\begin{thm}\label{pan4-2}\normalfont
If $n>1$, then
\begin{align*}
c_n(231,1423 : 231)=c_n(312,2314 : 312)&=
  \begin{cases}
    \frac{7}{3}\cdot4^{k-1}-\frac{1}{3} & \text{if } n=2k, \\
    \frac{14}{3}\cdot4^{k-1}-\frac{2}{3} & \text{if } n=2k+1.
  \end{cases}
\end{align*}
\end{thm}

\demo Note that $2314=comp(rev(1423))$. Thus, $c_n(231,1423 : 231)=c_n(312,2314 : 312)$ by Lemma \ref{pan2-1}. So it suffices to consider $c_n(312,2314 : 312)$. Next we discuss $c^{(i)}_n(312,2314 : 312)$ in different situations.

It follows from Theorem\ \ref{pan2-2} that $c^{(n)}_n(312,2314 : 312)=n-\lceil \frac{n}{2}\rceil$. For $i=1,2$, it is clear that $c^{(1)}_n(312,2314 : 312)=c_{n-1}(312,2314 : 312)$ and $c^{(2)}_n(312,2314 : 312)=c_{n-2}(312,2314 : 312)$. Let $\pi = \pi_1\pi_2\cdots \pi_n$ with $\pi_i=1$ and $2<i\leq n-1$. Note that $\pi\in C^{(i)}_n(312,2314 : 312)$ if and only if $\pi=\delta_i\oplus\mu$ where $\mu\in C_{n-i}(312,2314 : 312)$. Hence, $c^{(i)}_n(312,2314 : 312)=c_{n-i}(312,2314 : 312)$ for $2<i\leq n-1$, and so $$c_n(312,2314 : 312)=\sum_{i=1}^{n-1}c_{n-i}(312,2314 : 312)+n-\lceil \frac{n}{2}\rceil.$$
In addition, it is easy to derive that $$c_n(312,2314 : 312)=4c_{n-2}(312,2314 : 312)+3+2\lceil \frac{n-2}{2}\rceil-\lceil \frac{n-1}{2}\rceil-\lceil \frac{n}{2}\rceil.$$
Note that
\begin{align*}
2\lceil \frac{n-2}{2}\rceil-\lceil \frac{n-1}{2}\rceil-\lceil \frac{n}{2}\rceil&=
  \begin{cases}
    -2 & \text{if } n =2k, \\
    -1 & \text{if } n =2k+1.
  \end{cases}
\end{align*}
Note that $c_2(312,2314 : 312)=2$ and $c_3(312,2314 : 312)=4$. Thereby, we deduce that
\begin{align*}
c_n(312,2314 : 312)&=
  \begin{cases}
    \frac{7}{3}4^{k-1}-\frac{1}{3} & \text{if } n =2k \\
    \frac{14}{3}4^{k-1}-\frac{2}{3} & \text{if } n =2k+1
  \end{cases}
\end{align*}
The proof of this theorem is complete. \qed

\begin{thm}\label{pan4-3}\normalfont
If $n\geq1$, then
\begin{align*}
c_n(231,1243 : 231)=c_n(312,2134 : 312)&=
  \begin{cases}
    \frac{4k^3+3k^2-k}{6}+1 & \text{if } n =2k, \\
    \frac{4k^3+9k^2+5k}{6}+1 & \text{if } n =2k+1.
  \end{cases}
\end{align*}
\end{thm}
\demo By the same token, $c_n(231,1243 : 231)=c_n(312,2134 : 312)$ and it suffices to consider $c_n(312,2134 : 312)$. Next we discuss $c^{(i)}_n(312,2134 : 312)$ in different situations.

Applying Theorem\ \ref{pan2-2}, it follows that $c^{(n)}_n(312,2134 : 312)=n-\lceil \frac{n}{2}\rceil$. Moreover, it is clear that $c^{(1)}_n(312,2134 : 312)=c_{n-1}(312,2134 : 312)$, and $c^{(2)}_n(312,2134 : 312)=1$ because that $\pi\in C^{(2)}_n(312,2134 : 312)$ if and only if $\pi=21n(n-1)\cdots 3$. Let $\pi = \pi_1\pi_2\cdots \pi_n\in S_n$ with $\pi_i$ and $2<i\leq n-1$. Note that $\pi\in C^{(i)}_n(312,2134 : 312)$ if and only if $\pi=\mu\oplus\delta_{n-i}$ where $\mu\in C^{(i)}_i(312,2134 : 312)$. Hence, $c^{(i)}_n(312,2134 : 312)=c^{(i)}_{i}(312,2134 : 312)$ for $2<i\leq n-1$. Since $c^{(2)}_{2}(312,2134 : 312)=c^{(2)}_n(312,2134 : 312)=1$, we deduce that $$c_n(312,2134 : 312)=c_{n-1}(312,2134 : 312)+\sum_{i=2}^{n}c^{(i)}_{i}(312,2134 : 312),$$
namely,
$$c_n(312,2134 : 312)-c_{n-1}(312,2134 : 312)=\sum_{i=2}^{n}(i-\lceil \frac{i}{2}\rceil).$$
Since $$\sum_{i=2}^{n}(i-\lceil \frac{i}{2}\rceil)=\sum_{i=1}^{n}(i-\lceil \frac{i}{2}\rceil)-(1-\lceil \frac{1}{2}\rceil)=\sum_{i=1}^{n}(i-\lceil \frac{i}{2}\rceil)=\sum_{i=1}^{n}i-\sum_{i=1}^{n}\lceil \frac{i}{2}\rceil,$$
we have
$$c_n(312,2134 : 312)-c_{n-1}(312,2134 : 312)=\frac{n^2+n}{2}-\sum_{i=1}^{n}\lceil \frac{i}{2}\rceil,$$
$$c_{n-1}(312,2134 : 312)-c_{n-2}(312,2134 : 312)=\frac{(n-1)^2+(n-1)}{2}-\sum_{i=1}^{n-1}\lceil \frac{i}{2}\rceil,$$
$$\vdots$$
$$c_4(312,2134 : 312)-c_3(312,2134 : 312)=\frac{4^2+4}{2}-\sum_{i=1}^{4}\lceil \frac{i}{2}\rceil.$$
Therefore,
$$c_n(312,2134 : 312)-c_3(312,2134 : 312)=\sum_{i=4}^{n}\frac{i^2+i}{2}-\sum_{j=4}^{n}(\sum_{i=1}^{j}\lceil \frac{i}{2}\rceil).$$
Clearly,$$\sum_{i=4}^{n}\frac{i^2+i}{2}=\sum_{i=1}^{n}\frac{i^2+i}{2}-10=\frac{2n^3+6n^2+4n}{12}-10.$$
If $n=2k$, then by Lemma\ \ref{pan2-3} we deduce that
$$\sum_{j=4}^{n}(\sum_{i=1}^{j}\lceil \frac{i}{2}\rceil)=\sum_{s=2}^{k}(\sum_{i=1}^{2s}\lceil \frac{i}{2}\rceil)+\sum_{t=2}^{k-1}(\sum_{i=1}^{2t+1}\lceil \frac{i}{2}\rceil)=
\sum_{s=2}^{k}(s^2+s)+\sum_{t=2}^{k-1}(t^2+2t+1)=$$
$$\frac{4k^3+9k^2+5k}{6}-7.$$
Similarly, if $n=2k+1$, then
$$\sum_{j=4}^{n}(\sum_{i=1}^{j}\lceil \frac{i}{2}\rceil)=\sum_{j=4}^{n-1}(\sum_{i=1}^{j}\lceil \frac{i}{2}\rceil)+\sum_{i=1}^{n}\lceil \frac{i}{2}\rceil=\frac{4k^3+9k^2+5k}{6}-7+k^2+2k+1=$$
$$=\frac{4k^3+15k^2+17k}{6}-6.$$
Therefore,
\begin{align*}
c_n(312,2134 : 312)&=
  \begin{cases}
    \frac{4k^3+3k^2-k}{6}+1 & \text{if } n =2k, \\
    \frac{4k^3+9k^2+5k}{6}+1 & \text{if } n =2k+1.
  \end{cases}
\end{align*}
One easily checks that the formula also holds for $n=1,2,3$. The proof of this theorem is complete. \qed

\begin{thm}\label{pan4-4}\normalfont
If $n\geq1$, then
\begin{align*}
c_n(231,2143 : 231)=c_n(312,2143 : 312)&=
  \begin{cases}
     \frac{4k^3+3k^2-k}{6}+1 & \text{if } n =2k, \\
    \frac{4k^3+9k^2+5k}{6}+1 & \text{if } n =2k+1.
  \end{cases}
\end{align*}
\end{thm}
\demo By the same token, $c_n(231,2143 : 231)=c_n(312,2143 : 312)$ and it suffices to consider $c_n(312,2143 : 312)$. Next we discuss $c^{(i)}_n(312,2143 : 312)$ in different situations.

Applying Theorem\ \ref{pan2-2}, it follows that $c^{(n)}_n(312,2143 : 312)=n-\lceil \frac{n}{2}\rceil$. Moreover, it is clear that $c^{(1)}_n(312,2143 : 312)=c_{n-1}(312,2143 : 312)$, and $c^{(2)}_n(312,2143 : 312)=1$ because that $\pi\in C^{(2)}_n(312,2134 : 312)$ if and only if $\pi=2134\cdots n$. Let $\pi = \pi_1\pi_2\cdots \pi_n\in S_n$ with $\pi_i$ and $2<i\leq n-1$. Note that $\pi\in C^{(i)}_n(312,2143 : 312)$ if and only if $\pi=\mu\oplus id_{n-i}$ where $\mu\in C^{(i)}_i(312,2143 : 312)$. Hence, $c^{(i)}_n(312,2143 : 312)=c^{(i)}_{i}(312,2143 : 312)$ for $2<i\leq n-1$. Since $c^{(2)}_{2}(312,2143 : 312)=c^{(2)}_n(312,2143 : 312)=1$, we deduce that $$c_n(312,2143 : 312)=c_{n-1}(312,2143 : 312)+\sum_{i=2}^{n}c^{(i)}_{i}(312,2143 : 312).$$
Proceeding as in the proof of Theorem\ \ref{pan4-3}, we prove this theorem. \qed

In fact, we can deduce the rest cases for $\tau\in S_4$ by applying this method. For example, we know that $\pi\in C^{(i)}_n(312,1234 : 312)$ if and only if $\pi=\mu\oplus\tau$ where $\mu\in C^{(i)}_i(312,1234 : 312)$ and $\tau\in C_{n-i}(312,123 : 312)$, and therefore, we can obtain the recursive relationship by using Theorem\ \ref{pan3-1}.
However, the calculations are very complex and may not even be able to provide formulas, so we will not study them one by one.

%

\vskip0.4cm
\noindent{\bf Conflict of interest} The authors state no conflict of interest. In this paper, there is no experimental or computer calculated data, and there is only
logical proof.

\end{document}